\documentclass[11pt]{article}
\usepackage{amsfonts}

\usepackage{amssymb}
\usepackage{latexsym}
\usepackage{amscd}
\usepackage[centertags]{amsmath}
\usepackage{xypic}

\newtheorem{thm}{Theorem}

\newtheorem{lem}[thm]{Lemma}

\numberwithin{equation}{section}

\newcommand{\rmd}{\mathrm{d}}

\newcommand{\Rnum}{\mathbb{R}}

\newcommand{\Nnum}{\mathbb{N}}

\newcommand{\lam}{\lambda}

\newcommand{\vep}{\varepsilon}
\newcommand{\qed}{\hfill $\Box$}

\newcommand{\wt}{\widetilde}

\begin{document}
\title{Bounding Ornstein-Uhlenbeck Processes and Alikes}
\author{Jian-Sheng Xie\footnote{E-mail: jiansheng.xie@gmail.com}\\
{\footnotesize School of Mathematical Sciences, Fudan
University, Shanghai 200433, P. R. China}} 
\date{}
\maketitle
\begin{abstract}
In this note we consider SDEs of the type $\rmd X_t=[F (X_t) -A X_t] \rmd t +D \rmd W_t$ under
the assumptions that $A$'s eigenvalues are all of positive real parts and $F (\cdot)$ has
slower-than-linear growth rate. It is proved that $\displaystyle \varlimsup_{t \to
\infty} \frac{\|X_t\|}{\sqrt{\log t}} =\sqrt{2 \lam_1}$ almost surely with $\lam_1$ being 
the largest eigenvalue of the matrix $\displaystyle \Sigma :=\int_0^\infty e^{-s A} \cdot
(D \cdot D^T) \cdot e^{-s A^T} \rmd s$; the discarded measure-zero set can be chosen
independent of the initial values $X_0=x$.
\end{abstract}

\section{Introduction}\label{sec:1}
It's well known that, for a given one-dimensional stationary Ornstein-Uhlenbeck (OU for short) process
$X=\{X_t: t \geq 0\}$£¬there exist $\lam, \sigma>0, \mu \in \Rnum$ and a standard Brownian Motion (BM for
short) $B (\cdot)$ such that $X$ has the same distribution as $\{\sigma \cdot e^{-\frac{\lam t}{2}} \cdot
B (e^{\lam t}) + \mu: t \geq 0 \}$. Therefore the law of iterated logarithm for BM (see, e.g., \cite{RM99})
leads us to the conclusion $X_t= O (\sqrt{\log t})$ almost surely. In this note we investigate what bounds
can we achieve for higher dimensional OU processes $X=\{X_t: t \geq 0\}$ and alikes which may be modeled
by the following SDE (of dimension $d \geq 2$)
\begin{equation}\label{eq: SDE-0}
\rmd X_t=[F (X_t) -A X_t] \rmd t +D \, \rmd W_t,
\end{equation}
where $D$ is a constant $d$-by-$d$ matrix. And we always assume the following conditions:
\begin{itemize}
  \item[(\textbf{C1})] All the eigenvalues of the $d$-by-$d$ matrix $A$ have positive real parts;
  \item[(\textbf{C2})] $F (x)=o (\|x\|)$ (as $\|x\| \to \infty$) is a continuous $\Rnum^d$-valued function. 
Here $\|\cdot\|$ denotes the standard Euclidean norm.
\end{itemize}

Our main result can be stated as the following.
\begin{thm}\label{thm: main}
The solution to (\ref{eq: SDE-0}) always satisfies
\begin{equation}\label{eq: limsup}
\varlimsup_{t \to \infty} \frac{\|X_t\|}{\sqrt{\log t}} =\sqrt{2 \lam_1} \hbox{ almost surely, }
\end{equation}
where $\lam_1$ is the largest eigenvalue of the matrix $\displaystyle \Sigma :=\int_0^\infty e^{-s A} \cdot
(D \cdot D^T) \cdot e^{-s A^T} \rmd s$.
Here the discarded measure-zero set can be chosen independent of the initial values $X_0=x$.
\end{thm}
Such result seems to be new in literature as to our knowledge and deserves a publication somewhere.
The proof of the main theorem, based mainly on the well-known fact mentioned at the beginning of the 
introduction and on elemental linear algebra, is presented in Sect. \ref{sec:2} and Sect. \ref{sec:3}; 
the calculation of the precise limit value in (\ref{eq: limsup}) is based mainly on \cite{Xie14-2}, 
see Sect. \ref{sec:2}.
\section{OU Processes Case: $F =0$}\label{sec:2}
In this part, we consider the simpler case of $F=0$, i.e., the follow model
\begin{equation}\label{eq: SDE-1}
\rmd X_t= -A X_t \rmd t +D \, \rmd W_t.
\end{equation}
Clearly the solution satisfies the follow formula
\begin{equation}\label{eq: SDE-1'}
X_t=e^{-t A} X_0 +\int_0^t e^{-(t-s) A} D \, \rmd W_s.
\end{equation}
When $X_0 \sim N (0, \Sigma)$ with $\displaystyle \Sigma :=\int_0^\infty e^{-s A} \cdot
(D \cdot D^T) \cdot e^{-s A^T} \rmd s$, $\{X_t: t \geq 0\}$ is a stationary Markov process.

Throughout this section, we will use $B_{\cdot}$ in denoting one dimensional standard BM
and write $W_{\cdot}$ for higher dimensional standard BM.

As we have addressed in the introduction, any one dimensional stationary OU process is of
growth rate $O (\sqrt{\log t})$. This result can be restated as the following lemma, whose
proof is omitted.
\begin{lem}\label{lem: 1}
For any $\lam>0$, almost surely we have
$$
\int_0^t e^{-\lam (t-s)} \, \rmd B_s=O (\sqrt{\log t}).
$$
\end{lem}

Based on the above lemma, we would prove the following three lemmas.
\begin{lem}\label{lem: 2}
For $\lam>0$ and any $k \in \Nnum$, almost surely we have
$$
\int_0^t \frac{(t-s)^k}{k!} \cdot e^{-\lam (t-s)} \, \rmd B_s=O (\sqrt{\log t}).
$$
\end{lem}
\begin{lem}\label{lem: 3}
For $\lam>0, \mu \neq 0$, almost surely we have
\begin{eqnarray*}
\int_0^t e^{-\lam (t-s)} \cdot \cos \mu (t-s) \, \rmd B_s &=& O (\sqrt{\log t}),\\
\int_0^t e^{-\lam (t-s)} \cdot \sin \mu (t-s) \, \rmd B_s &=& O (\sqrt{\log t}).
\end{eqnarray*}
\end{lem}
\begin{lem}\label{lem: 4}
For $\lam>0, \mu \neq 0$ and any $k \in \Nnum$, almost surely we have
\begin{eqnarray*}
\int_0^t \frac{(t-s)^k}{k!} \cdot e^{-\lam (t-s)} \cdot \cos \mu (t-s) \, \rmd B_s &=& O (\sqrt{\log t}),\\
\int_0^t \frac{(t-s)^k}{k!} \cdot e^{-\lam (t-s)} \cdot \sin \mu (t-s) \, \rmd B_s &=& O (\sqrt{\log t}).
\end{eqnarray*}
\end{lem}

\noindent \textit{Proof of Lemma \ref{lem: 2}. \;}
Put $Y_t :=\int_0^t e^{-\lam (t-s)} \rmd B_s$ and
$$
L (t) :=\sup \limits_{u \in [0, t]} |Y_u|, \quad I^{(k)}_t :=\int_0^t \frac{(t-s)^k}{k!}
\cdot e^{-\lam (t-s)} \, \rmd B_s.
$$
Clearly
$$
I^{(1)}_t=\int_0^t e^{-\lam (t-u)} Y_u \rmd u, \quad I^{(k+1)}_t=\int_0^t e^{-\lam (t-u)}
I^{(k)}_u \rmd u, k \geq 1.
$$
And
$$
|I^{(1)}_t| \leq \int_0^t e^{-\lam (t-u)} |Y_u| \rmd u \leq \int_0^t e^{-\lam (t-u)} L (u) \rmd u \leq L (t)/\lam.
$$
Lemma \ref{lem: 1} tells us $L (t)=O (\sqrt{\log t})$. Hence $I^{(1)}_t=O (\sqrt{\log t})$.
Inductively $I^{(k)}_t=O (\sqrt{\log t})$ for all $k \in \Nnum$.
\qed

\noindent \textit{Proof of Lemma \ref{lem: 3}. \;}
For any $\theta \in \Rnum$, we write
$$
R_\theta := \left[
    \begin{array}{rl}
    \cos \theta \; & -\sin \theta\\
    \sin \theta \; & \cos \theta
    \end{array}
       \right].
$$
These are rotations which preserve the distance induced by the standard norm $\|\cdot\|$ on
$\Rnum^2$.

Consider the following diffusion process
$$
X_t :=\int_0^t e^{-\lam (t-s)} \cdot R_{-\mu (t-s)} \, \rmd W_s,
$$
where $W_{\cdot}=(W_{\cdot}^1, W_{\cdot}^2)^T$ is a 2-dimensional standard BM. Define
$$
\wt{W}_t :=\int_0^t R_{\mu s} \, \rmd W_s.
$$
It is easy to see that $\wt{W}_{\cdot}$ is still a 2-dimensional standard BM. And
$$
X_t=\int_0^t e^{-(t-s)} \cdot R_{-\mu t} \, \rmd \wt{W}_s.
$$
Now in view of Lemma \ref{lem: 1} it is clear that
$$
\|X_t\|=\|\int_0^t e^{-(t-s)} \rmd \wt{W}_s\|=O (\sqrt{\log t}).
$$
Thus
$$
\int_0^t e^{-\lam (t-s)} \cdot \Bigl[ \cos \mu (t-s) \, \rmd W_s^1 -\sin \mu (t-s)
\, \rmd W_s^2 \Bigr] =O (\sqrt{\log t}).
$$
Similarly,
$$
\int_0^t e^{-\lam (t-s)} \cdot \Bigl[ \cos \mu (t-s) \, \rmd W_s^1 +\sin \mu (t-s)
\, \rmd W_s^2 \Bigr] =O (\sqrt{\log t}).
$$
The lemma follows from the above equations.
\qed

\noindent \textit{Proof of Lemma \ref{lem: 4}. \;}
Now for any $k \geq 1$ (fixed), consider
$$
X_t :=\int_0^t \frac{(t-s)^k}{k!} \cdot e^{-\lam (t-s)} \cdot R_{-\mu (t-s)} \, \rmd W_s,
$$
where $R_\cdot$ is introduced in the proof of Lemma \ref{lem: 3}. It is easy
to see that
$$
\|X_t\|=\|\int_0^t \frac{(t-s)^k}{k!} \cdot e^{-\lam (t-s)} \, \rmd \wt{W}_s\|,
$$
where $\wt{W}_\cdot$ is also introduced in the proof of Lemma \ref{lem: 3}. Now Lemma
\ref{lem: 2} tells us $\|X_t\|=O (\sqrt{\log t})$ and the rest proof follows smoothly as in the
proof of Lemma \ref{lem: 3}.
\qed

From the above four lemmas we easily prove the bound $O (\sqrt{\log t})$ for the solutions to
(\ref{eq: SDE-1}) via exploiting the standard Jordan form of $A$ (and hence of $e^{-(t-s) A}$)
in formula (\ref{eq: SDE-1'}). The fact that the $\varlimsup$ in (\ref{eq: limsup}) is constant
almost surely follows from the ergodic property of the stationary OU process. 

Now we calculate the $\varlimsup$ in (\ref{eq: limsup}) 
explicitly via \cite{Xie14-2}: Without loss of generality, assume $\Sigma$ to be invertible. Take 
$V (x)=\frac{1}{2} x^T \Sigma^{-1} x$, the result in \cite{Xie14-2} tells us that
$\displaystyle \varlimsup_{t \to \infty} \frac{V (X_t)}{\sqrt{\log t}} \leq 1$ almost surely, 
which implies
$$
\varlimsup_{t \to \infty} \frac{\|X_t\|}{\sqrt{\log t}}=c \leq \sqrt{2 \lam_1} \hbox{ almost surely.}
$$
Now let $\alpha$ be a unit eigen-vector of $\Sigma$ corresponding to $\lam_1$. It is easy to see that
$Y:=\{Y_n:=\alpha^T X_n/{\sqrt{\lam_1}}: n \geq 0\}$ is a stationary Gaussian process with steady distribution 
$N (0, 1)$; This process inherits the exponential mixing property from $X_{\cdot}$. A standard result says 
that for i.i.d. standard normal random variables $\{Z_n: n \geq 0\}$, we always have $\displaystyle 
\varlimsup_{n \to \infty} \frac{|Z_n|}{\sqrt{2 \log n}}=1$ almost surely. With a tedious but routine effort 
(which we omit the details here), it is not hard to see that we still have $\displaystyle \varlimsup_{n \to 
\infty} \frac{|Y_n|}{\sqrt{2 \log n}} =1$ almost surely for the new process $Y$. Therefore
$$
c=\varlimsup_{t \to \infty} \frac{\|X_t\|}{\sqrt{\log t}} \geq \varlimsup_{n \to \infty} \frac{\|X_n\|}{\sqrt{\log n}}=\sqrt{2 \lam_1}.
$$
Hence $c=\sqrt{2 \lam_1}$. And (\ref{eq: limsup}) follows.
\section{General Case: $F \neq 0$}\label{sec:3}
Now we consider the general case with $F \neq 0$. As is known, the solution to
(\ref{eq: SDE-0}) satisfies
\begin{equation}\label{eq: solution2SDE-0}
X_t=e^{-t A} X_0 +\int_0^t e^{-(t-s) A} F (X_s) \rmd s +\int_0^t e^{-(t-s) A} D \, \rmd W_s.
\end{equation}
Define
$$
L (t)=\sup \limits_{u \in [0, t]} \|\int_0^u e^{-(u-s) A} D \, \rmd W_s\|.
$$
Clearly $L (t)=O (\sqrt{\log t})$.

Since $A$ satisfies condition (\textbf{C1}), there exist $\lam_0>0$ and $K>0$ such that
\begin{equation}\label{eq: 3.2}
\|e^{-t A}\| \leq K \cdot e^{-\lam_0 t}, \quad \forall t \geq 0.
\end{equation}
Fix an arbitrarily small $\vep>0$ (with $\vep<\frac{\lam_0}{K}$), by assumption (\textbf{C2})
there exists $C=C (\vep)>0$ such that
\begin{equation}\label{eq: 3.3}
\|F (x)\| \leq C +\vep \|x\|, \quad \forall x \in \Rnum^d.
\end{equation}
In view of (\ref{eq: solution2SDE-0}) we have
$$
\|X_t\| \leq K e^{-\lam_0 t} \|X_0\| +\int_0^t K  e^{-\lam_0 (t-s)} \Bigl(C +\vep \|X_s\|\Bigr) \rmd s
+L (t), \quad \forall t \geq 0.
$$
Define
$$
f (t) :=K \|X_0\|+L (t)+K C/\lam_0, \quad \vep_0 :=K \vep, \quad u (t) :=\|X_t\|.
$$
Then $u (\cdot)$ can be regarded as a continuous positive function (almost surely) which satisfies
the following inequality
\begin{equation}\label{eq: deterministic-0}
u (t) \leq f (t) +\vep_0 \int_0^t e^{-\lam_0 (t-s)} u (s) \rmd s, \quad \forall t \geq 0.
\end{equation}
Put $\displaystyle \varphi (t) :=\int_0^t e^{\lam_0 s} u (s) \rmd s$, we have
$$
\frac{\rmd \varphi}{\rmd t} (t) \leq f (t) \cdot e^{\lam_0 t} +\vep_0 \varphi (t), \quad t \geq 0
$$
which implies (where $\lam:=\lam_0-\vep_0$)
$$
\varphi (t) \leq e^{\vep_0 t} \cdot \int_0^t f (s) \cdot e^{\lam s} \rmd s, \quad t \geq 0.
$$
Thus, noting (\ref{eq: deterministic-0}) and the monotonicity of $f$, we have
\begin{eqnarray*}
u (t) &\leq& f (t)+\vep_0 e^{-\lam_0 t} \cdot \varphi (t) \leq f (t)+\vep_0 \int_0^t f (s) e^{-\lam (t-s)} \rmd s\\
&\leq& [1+\frac{\vep_0}{\lam_0-\vep_0}] \cdot f (t)=\frac{\lam_0}{\lam_0-\vep_0} \cdot f (t).
\end{eqnarray*}
This implies that the solution $X_t$ has almost the same growth rate as $\int_0^t e^{-(t-s) A} D \, \rmd W_s$.
Specifically we always have $\|X_t\|=O (\sqrt{\log t})$ almost surely. Clearly the limit value in
(\ref{eq: limsup}) is coincident with that of the OU case (i.e., the case $F=0$).

\noindent{\sl \textbf{Acknowledgements} \quad} {The author thanks Prof. Jiangang Ying for helpful discussions.
This work is partially supported by NSFC (No. 10701026 and No. 11271077) and the Laboratory of Mathematics for
Nonlinear Science, Fudan University.}






\begin{thebibliography}{00}

\bibitem{RM99} Revuz, Daniel; Yor, Marc:
{\it Continuous martingales and Brownian motion}.
Third edition. Grundlehren der Mathematischen Wissenschaften
[Fundamental Principles of Mathematical Sciences], 293.
Springer-Verlag, Berlin, 1999. xiv+602 pp. ISBN: 3-540-64325-7

\bibitem{Xie14-2} Xie, Jian-Sheng:
{\em Bounding the Solutions to Some SDEs via Ergodic Theory},
arXiv: 1407.2728.
%
%
%
%
%
%
%
%
%
%
%
%
%
%
%
%
%
%
%
%
%
%
%
%
%
%
%
%
%
%
%
%
%
%

\end{thebibliography}



\end{document}